\documentclass[11pt, reqno]{amsart}
\usepackage[utf8]{inputenc}
\usepackage{amsthm}
\usepackage{amsmath}
\usepackage{amssymb}
\usepackage{comment}
\usepackage[margin=1.5in]{geometry}
\usepackage{comment}
\usepackage[dvipsnames]{xcolor}
\usepackage{xcolor}

\usepackage[colorlinks=true, pdfstartview=FitV, linkcolor=blue, citecolor=blue, urlcolor=blue]{hyperref}

\providecommand{\C}{\mathbb C}

\providecommand{\N}{\mathbb N}

\providecommand{\R}{\mathbb R}
\providecommand{\Z}{\mathbb Z}

\newcommand{\E} {{\mathcal E}}

\newcommand{\norm}[1]{\left|\left| #1 \right|\right|}

\newtheorem{theorem}{Theorem}[section]

\newtheorem{lemma}[theorem]{Lemma}
\newtheorem{prop}[theorem]{Proposition}

\numberwithin{equation}{section}

\title{Sharpness of the Kadec'-type theorem for the union of intervals}
\author{Thibaud Alemany}
\address{Thibaud Alemany,  University of Central Florida,
	Department of Mathematics,
	Orlando, FL 32816, USA}
\email{thibaud.alemany@ucf.edu}

\author{Oleg Asipchuk}
\address{Oleg Asipchuk,  University of Cincinnati,
	Department of Mathematical Sciences,
	Cincinnati, OH 45221, USA}
\email{asipchah@ucmail.uc.edu}


\begin{document}
	
\begin{abstract}
	In this article, we establish the sharpness of a Kadec'-type theorem for finite unions of intervals admitting the integers as an orthonormal spectrum by constructing an explicit counterexample.
\end{abstract}

\maketitle

\section{Introduction}

We recall that an exponential basis on a domain $D\subset \R^d$ is a Riesz basis for $L^2(D)$ in the form of  \(\mathcal E(\Lambda)=\{e^{2\pi i\lambda\cdot x}:\lambda\in\Lambda\}\), where \(\Lambda\subset\mathbb R^d\) is a countable set of frequencies.   An important example of an exponential basis is the Fourier basis $ \E(\Z)$  for $L^2(0,1)$.

The stability of exponential bases under perturbations goes back to the classical work of Paley and Wiener, who established a general perturbation theorem for bases in Hilbert spaces, with applications to nonharmonic Fourier series \cite{PaleyWiener1934}. In 1964, M.~I.~Kadec' established sharp bounds on perturbations of orthogonal exponential bases for $L^2(0,1)$.
\begin{theorem}[\cite{Kadec1964}]\label{T-Kadec}
	Let $\Lambda =  \{\lambda_n \}_{n\in\Z}  $ be a sequence in $\R$ for which $$ | \lambda_n - n| \leq L < \frac{1}{4}$$ whenever $n\in\Z.$ Then $ \E(\Lambda)$ is an exponential basis for $L^2(0,1)$ with frame constants $A=(\cos(\pi L)-\sin( \pi L))^2$ and $B= (2-\cos(\pi L)+\sin( \pi L))^2$.   The constant $\frac{1}{4}$ cannot be replaced by any larger constant. 
\end{theorem}

This result was subsequently extended in several directions. Avdonin \cite{Avdonin1974} replaced the pointwise perturbation condition by a weaker condition controlling averages of the frequency perturbations over finite blocks. The frame bounds for this generalization were established in \cite{Alemany2025}.  More general characterizations of exponential Riesz bases were later obtained by Pavlov \cite{Pavlov1979} and by Khrushchev, Nikolskii, and Pavlov~\cite{Khrushchev1981}.

The stability problem was also extended beyond the one-dimensional interval. In \cite{Sun1999OnK1}, multidimensional and Fourier-frame analogues of Kadec' theorem were developed, while \cite{Balan1997} established related perturbation results for Fourier frames and more general frame systems. More recently, Kadec'-type stability principles have also been developed for unions of intervals \cite{asipchuk_drezels_2023}, Gabor systems \cite{DRV2026}, polynomial sampling systems \cite{MarzoSeip2009}, Fourier matrices \cite{Asipchuk2025}, and sampled group orbits \cite{KrishtalMiller2025}.

Let $M\ge 2$, $$0=a_1<a_2<...<a_{M}<a_{M+1}=1,$$ 
and $$-k=l_1<l_2<...<l_M=k-1,$$ where $k\in\N$ and $l_{j}\in \Z$. We define a finite union of the intervals
\begin{equation}\label{E-J}
    J=\bigcup_{j=1}^M (l_j+a_j,l_j+a_{j+1}).
\end{equation}
Note that $\E(\Z)$ is an orthonormal exponential basis for $L^2(J)$; see Proposition~\ref{P-ON}. Therefore, a slight modification of the proof of Theorem \ref{T-Kadec}, or alternatively the approach in \cite{KrishtalMiller2025}, yields the following Kadec'-type theorem for $J$.
\begin{theorem}\label{T-KI}
    Let $J$ be as in \eqref{E-J} and $\Lambda =  \{\lambda_n \}_{n\in\Z}  $ be a sequence in $\R$ for which $$ | \lambda_n - n| \leq L < \frac{1}{8k}$$ whenever $n\in\Z.$ Then $ \E(\Lambda)$ is an exponential basis for $L^2(J)$ with frame constants $A=(\cos(2kL\pi)-\sin( 2kL\pi))^2$ and $B= (2-\cos(2kL\pi)+\sin( 2kL\pi))^2$.
\end{theorem}
For the convenience of the reader, we provide a proof of Theorem~\ref{T-KI} in the appendix.

In this manuscript, we show that the bound $\tfrac{1}{8k}$ is optimal by constructing a counterexample. Our main result is the following.
\begin{theorem}\label{T-main}
    The bound $\tfrac{1}{8k}$ in Theorem \ref{T-KI} is sharp.
\end{theorem}
\medskip
The paper is organized as follows. In Section 2, we recall the necessary preliminaries and establish several auxiliary lemmas. In Section 3, we prove our main results. In the appendix, we provide a proof of Theorem~\ref{T-KI}. 

\medskip
\noindent
{\it AI Declaration}. During the preparation of this manuscript, the authors used ChatGPT, developed by OpenAI,  as an assistive tool for editing and typesetting the manuscript and for checking the correctness of mathematical arguments and proofs. All mathematical results, however, were obtained independently by the authors without the assistance of artificial intelligence.

\section{Preliminaries}
Let  $H$ be a  separable Hilbert space with inner product $\langle \ , \ \rangle $  and norm $||w ||=\sqrt{\langle w, w    \rangle} $. We will mostly work with $L^2(D),$ where $D\subset\R^d.$ Thus, the norm is given by $||f||_2^2=\int_D |f(x)|^2 dx.$   
  
A sequence of vectors ${\mathcal V}= \{v_j\}_{j\in\Z }$  in   $H$ is a
{\it Riesz basis} if 
there exist  constants $A, \ B>0$    such that, for any $w\in H$ and  for all finite sequences $ \{a_j\}_{j\in J}\subset\C $, the following inequalities hold. 
\begin{align}\label{e2-RS}
 	A  \sum_{j\in J}   |a_j|^2   \leq  &\Big\Vert \sum_{j\in J}  a_j  v_j \Big\Vert^2  \leq B \sum_{j\in J} |a_j|^2.
\\\label{e2-frame}
 A||w||^2\leq  &\sum_{j\in\Z} |\langle  w, v_j\rangle |^2\leq B ||w||^2.
 \end{align}
The constants $A$ and $B$ are called the {\it frame bounds}. If condition \eqref{e2-RS} holds, we call ${\mathcal V}$ a Riesz sequence, while if condition \eqref{e2-frame} holds, we call ${\mathcal V}$ a frame. If both \eqref{e2-RS} and \eqref{e2-frame} hold with $A=B=1$, then ${\mathcal V}$ is an orthonormal basis.

\begin{prop}\label{P-ON}
    $\E(\mathbb{Z})$ is an orthonormal exponential basis for $L^2(J)$, where $J$ is defined as in \eqref{E-J}.
\end{prop}

\begin{proof}   
Define $U:L^2(J)\to L^2(0,1)$ by
\begin{align*}
    (Uf)(t)=f(t+l_j),\quad t\in(a_j,a_{j+1}).
\end{align*}
Since the intervals $(a_j,a_{j+1})$ partition $(0,1)$ up to
finitely many endpoints, translation invariance of Lebesgue
measure gives
\begin{align*}
\|Uf\|_{L^2(0,1)}^2
=\sum_{j=1}^M\int_{a_j}^{a_{j+1}}|f(t+l_j)|^2\,dt
=\|f\|_{L^2(J)}^2.
\end{align*}
In particular, since $U$ is an isometry, it is one-to-one. Moreover, $U$ is onto: for $g\in L^2(0,1)$, define
$f(x)=g(x-l_j)$ on $(l_j+a_j,l_j+a_{j+1})$.
Thus $U$ is unitary. For every $n\in\mathbb Z$,  $l_j\in \Z$ implies
\begin{align*}
U\bigl(e^{2\pi i n\,\cdot}\bigr)(t)
=e^{2\pi i n(t+l_j)}
=e^{2\pi i nt}
\quad\text{for a.e. }t\in(0,1).
\end{align*}
Hence $U$ maps $\E(\mathbb Z)$ onto the standard Fourier
orthonormal basis of $L^2(0,1)$.
\end{proof}

Our goal is to show that a certain exponential system does not form a Riesz sequence on a union of intervals. To this end, we need estimates for particular double sums, which we establish in the following lemmas.
\begin{lemma}\label{Sharp Hilbert inequality}
For $N\in\mathbb{N}$, define 
    \begin{equation}\label{E-HN}
     H_N=\sum_{k=1}^N\frac{1}{k}.    
    \end{equation}
Then
\begin{equation}
\lim_{N\to\infty}
\frac{1}{H_N}
\sum_{n,m=1}^N\frac{1}{\sqrt{nm}(n+m)}
=\pi.
\end{equation}
\end{lemma}

\begin{proof}
Write
\begin{align*}
    S_N=\sum_{n,m=1}^N\frac{1}{\sqrt{nm}(n+m)}.
\end{align*}
For fixed $n$, the function
$y\mapsto y^{-1/2}(n+y)^{-1}$ is decreasing. Hence
\begin{align*}
    \sum_{m=1}^N\frac{1}{\sqrt m(n+m)}
\le
\int_0^\infty\frac{dy}{\sqrt y(n+y)}
=\frac{\pi}{\sqrt n},
\end{align*}
where the last equality follows by substituting $y=nu^2$.
Multiplying by $n^{-1/2}$ and summing gives $S_N\le \pi H_N$.

For the lower bound, the function
$(x,y)\mapsto(\sqrt{xy}(x+y))^{-1}$
is decreasing in each variable. Therefore, with $M=N+1$,
\begin{align*}
    S_N\ge
\int_1^M\int_1^M\frac{dy\,dx}{\sqrt{xy}(x+y)}.
\end{align*}
Integrating first in $y$, we obtain
\begin{align*}
S_N
&\ge
2\int_1^M
\left(
\arctan\sqrt{\frac Mx}
-\arctan\frac1{\sqrt x}
\right)\frac{dx}{x}\\
&=
\pi\log M
-2\int_1^M
\left(
\arctan\sqrt{\frac xM}
+\arctan\frac1{\sqrt x}
\right)\frac{dx}{x}.
\end{align*}
Here we used $\arctan u+\arctan(1/u)=\pi/2$ for $u>0$.
Since $\arctan u\le u$,
\begin{align*}
    \int_1^M\arctan\sqrt{\frac xM}\frac{dx}{x}
\le \frac1{\sqrt{M}}\int_1^M x^{-1/2}\,dx\le2,
\end{align*}
and
\begin{align*}
    \int_1^M\arctan\frac1{\sqrt{x}}\frac{dx}{x}
\le\int_1^M x^{-3/2}\,dx\le2.
\end{align*}
Consequently,
\begin{align*}
    \pi\log(N+1)-8\le S_N\le\pi H_N.
\end{align*}
Dividing by $H_N$ and using
$H_N=\log N+O(1)$ proves the limit.
\end{proof}

\begin{lemma}\label{L-sum-estimation}
    Let \(l>0\) be fixed. Then
    \begin{equation}\label{E-sum-estimation}
        \lim_{N\rightarrow \infty} \frac{1}{H_N}\sum_{n=1}^N\sum_{m=1}^N \frac{1}{\sqrt{nm}(n+m+l)}=\pi,       
    \end{equation}
\end{lemma}
\begin{proof}
   First, we observe that 
   \begin{align*}
       0&\le \frac{1}{H_N}\sum_{n=1}^N\sum_{m=1}^N \frac{1}{\sqrt{nm}(n+m)} -  \frac{1}{H_N}\sum_{n=1}^N\sum_{m=1}^N \frac{1}{\sqrt{nm}(n+m+l)}\\
       &\le\frac{l}{H_N} \sum_{n=1}^N\sum_{m=1}^N  \frac{1}{\sqrt{nm}(n+m)^2}\rightarrow 0,
   \end{align*}
   when $N\rightarrow\infty$, so it is enough to consider $l=0$. Now, the conclusion follows from Lemma \ref{Sharp Hilbert inequality}.  
\end{proof}

\begin{lemma}\label{L-cos}
Let $a,l\in(0,1)$. Then
\[
T_N=\sum_{n,m=1}^N
\frac{e^{2\pi ia(n+m)}}{\sqrt{nm}(n+m+l)}
\]
is uniformly bounded in $N$. More precisely,
\[
\sup_{N\ge1}|T_N|
\le
\frac{1}{|\sin(\pi a)|}
\sum_{m=1}^{\infty}\frac{1}{m^{3/2}}
<\infty.
\]
\end{lemma}

\begin{proof}
Since $a\in(0,1)$, the geometric series formula gives
\begin{align*}
    \left|\sum_{n=1}^k e^{2\pi ian}\right|
    \le \frac{1}{|\sin(\pi a)|},
    \qquad k\ge1.
\end{align*}
For each fixed $m\ge1$, the sequence
\begin{align*}
    b_n=\frac{1}{\sqrt{n}(n+m+l)}
\end{align*}
is positive and decreasing. Hence, by Abel summation,
\begin{align*}
\left|\sum_{n=1}^N e^{2\pi ian}b_n\right|
&\le \frac{1}{|\sin(\pi a)|}
\left(b_N+\sum_{n=1}^{N-1}(b_n-b_{n+1})\right)\\
&=\frac{b_1}{|\sin(\pi a)|}
=\frac{1}{(m+1+l)|\sin(\pi a)|}.
\end{align*}
Consequently,
\begin{align*}
|T_N|
&\le \sum_{m=1}^N\frac{1}{\sqrt m}
\left|\sum_{n=1}^N
\frac{e^{2\pi ian}}{\sqrt n(n+m+l)}\right|\\
&\le \frac{1}{|\sin(\pi a)|}
\sum_{m=1}^N\frac{1}{\sqrt m(m+1+l)}\\
&\le \frac{1}{|\sin(\pi a)|}
\sum_{m=1}^{\infty}\frac{1}{m^{3/2}}
<\infty.
\end{align*}
The bound is independent of $N$, which completes the proof.
\end{proof}

\section{Proof of the main theorem}
\begin{proof}[Proof of Theorem \ref{T-main}]
Let $J$ be defined as in \eqref{E-J}. Our counterexample is motivated by the classical counterexample demonstrating the sharpness of Kadec' $1/4-$theorem \cite{Levinson1936,RedhefferYoung1983}. We let  
\[
\lambda_n=
\begin{cases}
n+\frac{1}{8k}, & n>0,\\[4pt]
0, & n=0,\\[4pt]
n-\frac{1}{8k}, & n<0.
\end{cases}
\]
Next, for $N\in\mathbb{N}$,  we introduce a sequence $\mathcal{C}_N=\{c_n\}_{n}$ as
\begin{equation*}
    c_n=
\begin{cases}
\frac{1}{\sqrt{2nH_N}}, & 0<n\le N,\\[4pt]
0, & n=0\,\text{or}\,|n|>N\\[4pt]
-\frac{1}{\sqrt{-2nH_N }}, & -N\le n < 0,
\end{cases}
\end{equation*}
where $H_N$ is defined as in \eqref{E-HN}. In this case,
\[
\sum_n |c_n|^2
=
2\sum_{n=1}^N\frac{1}{2nH_N}
=
\frac{1}{H_N}\sum_{n=1}^N\frac1n
=
1.
\]
Next, we compute the quantity in \eqref{e2-RS} for $\mathcal{C}_N$ and $\E(\Lambda)$

\begin{align*}
    I=\left\|
\sum_n c_n e^{2\pi i\lambda_n x}
\right\|_J^2
&=
\left\|
\sum_{n=1}^N
\frac{e^{2\pi i(n+\frac{1}{8k})x}}
{\sqrt{2nH_N}}
-
\sum_{n=1}^N
\frac{e^{-2\pi i(n+\frac{1}{8k})x}}
{\sqrt{2nH_N}}
\right\|_J^2
\\
&=
2\left\|
\sum_{n=1}^N
\frac{\sin\!\left(2\pi(n+\frac{1}{8k})x\right)}
{\sqrt{nH_N }}
\right\|_J^2
\\
\end{align*}
Expanding the norm, we get
\begin{align*}
I&=
\frac{2}{H_N}
\sum_{n=1}^N\sum_{m=1}^N
\frac{1}
{\sqrt{nm}} \int_J
\sin\!\left(2\pi(n+\tfrac{1}{8k})x\right)
\sin\!\left(2\pi(m+\tfrac{1}{8k})x\right)\,dx
\\
&=\frac{1}{H_N}
\sum_{n=1}^N\sum_{m=1}^N
\frac{1}
{\sqrt{nm}} \int_J
[\cos\!\left(2\pi(n-m)x\right)-\cos\!\left(2\pi(n+m+\tfrac{1}{4k})x\right)]\,dx
\end{align*}
Next, we compute the integrals over each of the intervals comprising $J$. We first consider the case $m=n$. Then,
\begin{align*}
    \int_{l_j+a_j}^{l_j+a_{j+1}}
\cos\!\left(2\pi(n-m)x\right)\,dx&=\int_{l_j+a_j}^{l_j+a_{j+1}}
\,dx=a_{j+1}-a_j.
\end{align*}
Thus,
\begin{align}\label{E-nm}
    \int_J
\cos\!\left(2\pi(n-m)x\right)\,dx&=1.
\end{align}
Next, let $n\ne m$. In this case we have 
\begin{align*}
    \int_{l_j+a_j}^{l_j+a_{j+1}}
\cos\!\left(2\pi(n-m)x\right)\,dx&=\frac{\sin{(2\pi a_{j+1}(n-m))}-\sin{(2\pi a_{j}(n-m))}}{2\pi (n-m)}.
\end{align*}
Thus,
\begin{align}\label{E-nNotm}
    \int_J
\cos\!\left(2\pi(n-m)x\right)\,dx&=0.
\end{align}
Next, we directly compute the remaining integrals.
\begin{align*}
    \int_{l_j+a_j}^{l_j+a_{j+1}}
\cos\!\left(2\pi(n+m+\tfrac{1}{4k})x\right)\,dx&=\frac{\sin{(2\pi a_{j+1}(n+m+\tfrac{1}{4k})+\tfrac{l_j\pi }{4k})}}{2\pi (n+m+\tfrac{1}{4k})}\\
&-\frac{\sin{(2\pi a_{j}(n+m+\tfrac{1}{4k})+\tfrac{l_j\pi}{4k})}}{2\pi (n+m+\tfrac{1}{4k})}\
\end{align*}
Taking into account that $a_1=0$, $a_{M+1}=1$, $l_1=-k$, and $l_M=k-1$, we obtain
\begin{align*}
    \int_J
\cos\!\left(2\pi(n+m+\tfrac{1}{4k})x\right)\,dx&=\frac{1}{\pi (n+m+\tfrac{1}{4k})}\\
&-\sum_{j=2}^M \frac{\cos{(2\pi a_{j}(n+m+\tfrac{1}{4k})+\tfrac{(l_j+l_{j-1})\pi}{4k}) }\cdot\sin{(\tfrac{(l_j-l_{j-1})\pi}{4k})}}{\pi (n+m+\tfrac{1}{4k})}
\end{align*}
Combining this with \eqref{E-nm} and \eqref{E-nNotm}, we obtain
\begin{align*}
    I&=1- \frac{1}{\pi H_N}\sum_{n=1}^N\sum_{m=1}^N
\frac{1}{\sqrt{nm}(n+m+\tfrac{1}{4k})}\\
&+\sum_{j=2}^M\frac{\sin{(\tfrac{(l_j-l_{j-1})\pi}{4k})}}{\pi H_N}\sum_{n=1}^N\sum_{m=1}^N
\frac{\cos{(2\pi a_{j}(n+m+\tfrac{1}{4k})+\tfrac{(l_j+l_{j-1})\pi}{4k}) }}{\sqrt{nm}(n+m+\tfrac{1}{4k})}\\
&=1-\frac{I_1}{\pi}+\sum_{j=2}^M\frac{I_j}{\pi H_N}
\end{align*}
By Lemma \ref{L-sum-estimation}, $I_1\rightarrow \pi$, when $N\rightarrow \infty$. By Lemma \ref{L-cos}, the \(I_j\) are bounded by constants that do not depend on \(N\), so 
$\frac{\sum_{j=2}^{M} I_j}{\pi H_N}\rightarrow 0$, when $N\rightarrow \infty$. Therefore, $I\rightarrow 0$, when $N\rightarrow \infty$. Hence, $\E(\Lambda)$ is not a Riesz sequence and, consequently, cannot form a Riesz basis.
\end{proof}

  {\bf Remark.} The method described in the proof of Theorem \ref{T-main} can also be used to reprove the sharpness of Theorem \ref{T-Kadec}. Indeed, we let  
\[
\lambda_n=
\begin{cases}
n+\frac{1}{4}, & n>0,\\[4pt]
0, & n=0,\\[4pt]
n-\frac{1}{4}, & n<0.
\end{cases}
\]  
Next, for each $N\in\mathbb{N}$, we consider the sequence $\mathcal{C}_N$ defined as
\begin{equation*}
    c_n=
\begin{cases}
\frac{(-1)^n}{\sqrt{2nH_N}}, & 0<n\le N,\\[4pt]
0, & n=0\,\text{or}\,|n|>N\\[4pt]
-\frac{(-1)^n}{\sqrt{-2nH_N }}, & -N\le n < 0,
\end{cases}
\end{equation*}
So, we compute the quantity in \eqref{e2-RS} for $\mathcal{C}_N$ and $\E(\Lambda)$
\begin{align*}
    I=\left\|
\sum_n c_n e^{2\pi i\lambda_n x}
\right\|_J^2
&=
\left\|
\sum_{n=1}^N
\frac{(-1)^n e^{2\pi i(n+\frac{1}{4})x}}
{\sqrt{2nH_N}}
-
\sum_{n=1}^N
\frac{(-1)^n e^{-2\pi i(n+\frac{1}{4})x}}
{\sqrt{2nH_N}}
\right\|_J^2
\\
&=
2\left\|
\sum_{n=1}^N
\frac{(-1)^n\sin\!\left(2\pi(n+\frac{1}{4})x\right)}
{\sqrt{nH_N }}
\right\|_J^2,
\\
\end{align*}
where $J$ is $[-1/2,1/2]$. Expanding the norm, we get
\begin{align*}
I&=
\frac{2}{H_N}
\sum_{n=1}^N\sum_{m=1}^N
\frac{(-1)^{n+m}}
{\sqrt{nm}} \int_J
\sin\!\left(2\pi(n+\tfrac{1}{4})x\right)
\sin\!\left(2\pi(m+\tfrac{1}{4})x\right)\,dx
\\
&=\frac{1}{H_N}
\sum_{n=1}^N\sum_{m=1}^N
\frac{(-1)^{n+m}}
{\sqrt{nm}} \int_J
[\cos\!\left(2\pi(n-m)x\right)-\cos\!\left(2\pi(n+m+\tfrac{1}{2})x\right)]\,dx\\
&=1- \frac{1}{\pi H_N}\sum_{n=1}^N\sum_{m=1}^N
\frac{1}{\sqrt{nm}(n+m+\tfrac{1}{2})}=1-\frac{I_1}{\pi}.
\end{align*}
By Lemma \ref{L-sum-estimation}, $I_1\rightarrow \pi$, when $N\rightarrow \infty$. Therefore, $I\rightarrow 0$, when $N\rightarrow \infty$. Hence, $\E(\Lambda)$ is not a Riesz sequence and, consequently, cannot form a Riesz basis.

\appendix
\section{Proof of Theorem~\ref{T-KI}}
We adapt Kadec' original proof, which can also be found in Young's book \cite[Theorem 14]{Young2001}. The proof relies on the following theorem; see \cite[Theorem 10]{Young2001}.

\begin{theorem}[Paley--Wiener]\label{T-PW}
	Let $\{x_n\}_{n\in\N}$ and $\{y_n\}_{n\in\N}$ be sequences in a Banach space $X$. Let  $\lambda$ be a real number $(0<\lambda<1)$ such that 
	\begin{equation*}
		\norm{\sum_n c_n(x_n-y_n)}\leq \lambda \norm{\sum_n c_n x_n} 
	\end{equation*}
	holds for every finite sequence of scalars $\{c_n\}\subset \C.$ Then if $\{x_n\}$ is a basis so is $\{y_n\}.$ Moreover, if $\{x_n\}$ has Riesz constants $A$ and $B,$ then 
	$$
	(1-\lambda)^2 A  \sum   |c_n|^2   \leq  \left\Vert \sum  c_n  y_n \right\Vert^2  \leq (1+\lambda)^2 B \sum |c_n|^2. 
	$$
\end{theorem}

\begin{proof}[Proof of Theorem \ref{T-KI}]
    It is to be shown that
\[
\left\|
\sum c_n\left(e^{2\pi int}-e^{2\pi i\lambda_n t}\right)
\right\|_J
\leq \lambda<1
\]
whenever
\(
\sum |c_n|^2\leq1.
\)
Write
\[
e^{2 \pi int}-e^{2\pi i\lambda_n t}
=
e^{2\pi int}\left(1-e^{2\pi i\delta_n t}\right),
\]
where
\(
\delta_n=\lambda_n-n,
\) $n\in\Z$. The trick is to expand the function
$1-e^{2\pi ibt},$ $-k\leq t\leq k$ in a Fourier series relative to the complete orthogonal system
\[
\{1,\cos {(\tfrac{\pi lt}{k})},\sin((l-\tfrac12)\frac{\pi t}{k})\}_{l=1}^{\infty}
\]
and then exploit the fact that $|\lambda_n-n|$ is not too large.

Simple calculations show that when $\delta$ is real,
\begin{align*}
1-e^{2 \pi i\delta t}
&=
\left(1-\frac{\sin(2\pi k\delta}{2\pi k\delta}\right)
+
\sum_{l=1}^{\infty}
\frac{(-1)^l 4k\delta\sin(2\pi k\delta)}
     {\pi(l^2-4k^2\delta^2)}
\cos \tfrac{\pi l t}{k}\\
&+i\sum_{l=1}^{\infty}
\frac{(-1)^l4 k\delta\cos(2\pi k\delta)}
     {\pi\left((l-1/2)^2-4 k^2\delta^2\right)}
\sin((l-\tfrac12)\tfrac{\pi t}{k}).
\end{align*}
Let $\{c_n\}$ be an arbitrary finite sequence of scalars such that
\(
\sum |c_n|^2\leq1.
\)
By interchanging the order of summation and then using the triangle
inequality we see that
\[
\left\|
\sum c_n\left(e^{2 \pi int}-e^{2 \pi i\lambda_n t}\right)
\right\|
\leq A+B+C,
\]
where
\[
A=
\left\|
\sum_n
\left(
1-\frac{\sin(2 \pi k\delta_n)}{2\pi k\delta_n}
\right)c_ne^{2 \pi int}
\right\|,
\]
\[
B=
\sum_{l=1}^{\infty}
\left\|
\cos (\tfrac{\pi lt}{k})
\sum_n
\frac{(-1)^l4 k \delta_n\sin(2\pi k\delta_n)}
     {\pi(l^2-4 k^2\delta_n^2)}
c_ne^{2\pi i nt}
\right\|,
\]
and
\[
C=
\sum_{l=1}^{\infty}
\left\|
\sin((l-\tfrac12)\tfrac{\pi t}{k})
\sum_n
\frac{(-1)^l4 k\delta_n\cos(2\pi k\delta_n)}
     {\pi\left((l-1/2)^2-4 k^2\delta_n^2\right)}
c_ne^{2 \pi i nt}
\right\|.
\]
Obvious estimates show that
\[
A\leq
\left(1-\frac{\sin(2 \pi kL)}{2\pi kL}\right),\quad B\leq
\sum_{l=1}^{\infty}
\frac{4k L\sin(2\pi kL)}{\pi(l^2-4 k^2L^2)},
\]
and
\[
C\leq
\sum_{l=1}^{\infty}
\frac{4kL\cos(2\pi kL)}
{\pi\left((l-\tfrac12)^2-4k^2L^2\right)}.
\]
But the series
\[
\sum_{l=1}^{\infty}
\frac{4kL}{\pi(l^2-4k^2L^2)}\quad \text{and} \quad \sum_{l=1}^{\infty}
\frac{4kL}
{\pi\left((l-\tfrac12)^2-4k^2L^2\right)}
\]
are the partial fraction expansions of the functions
\[
\frac{1}{2 \pi kL}-\cot(2\pi kL)
\qquad\text{and}\qquad
\tan{(2\pi kL)},
\]
respectively, so that
\[
\left\|
\sum c_n\left(e^{2\pi int}-e^{2\pi i\lambda_n t}\right)
\right\|_J
\leq
\lambda
=
1-\cos{(2\pi kL)}+\sin{(2\pi kL)}.
\]
Therefore,  $L<\frac{1}{8k}$ implies $\lambda<1$.
\end{proof}

\bibliographystyle{plain}

\end{document}